\magnification 1250
\pretolerance=500 \tolerance=1000  \brokenpenalty=5000
\mathcode`A="7041 \mathcode`B="7042 \mathcode`C="7043
\mathcode`D="7044 \mathcode`E="7045 \mathcode`F="7046
\mathcode`G="7047 \mathcode`H="7048 \mathcode`I="7049
\mathcode`J="704A \mathcode`K="704B \mathcode`L="704C
\mathcode`M="704D \mathcode`N="704E \mathcode`O="704F
\mathcode`P="7050 \mathcode`Q="7051 \mathcode`R="7052
\mathcode`S="7053 \mathcode`T="7054 \mathcode`U="7055
\mathcode`V="7056 \mathcode`W="7057 \mathcode`X="7058
\mathcode`Y="7059 \mathcode`Z="705A
\def\spacedmath#1{\def\packedmath##1${\bgroup\mathsurround =0pt##1\egroup$}
\mathsurround#1
\everymath={\packedmath}\everydisplay={\mathsurround=0pt}}
\def\nospacedmath{\mathsurround=0pt
\everymath={}\everydisplay={} } \spacedmath{2pt}
\def\qfl#1{\buildrel {#1}\over {\longrightarrow}}
\def\hfl#1#2{\normalbaselines{\baselineskip=0truept
\lineskip=10truept\lineskiplimit=1truept}\nospacedmath\smash{\mathop{\hbox to
12truemm{\rightarrowfill}}\limits^{\scriptstyle#1}_{\scriptstyle#2}}}
\def\diagram#1{\def\normalbaselines{\baselineskip=0truept
\lineskip=10truept\lineskiplimit=1truept}   \matrix{#1}}
\def\vfl#1#2{\llap{$\scriptstyle#1$}\left\downarrow\vbox to
6truemm{}\right.\rlap{$\scriptstyle#2$}}
\def\mono{\lhook\joinrel\mathrel{\longrightarrow}}
\def\iso{\vbox{\hbox to .8cm{\hfill{$\scriptstyle\sim$}\hfill}
\nointerlineskip\hbox to .8cm{{\hfill$\longrightarrow $\hfill}} }}

\def\sdir_#1^#2{\mathrel{\mathop{\kern0pt\oplus}\limits_{#1}^{#2}}}
\def\pprod_#1^#2{\raise
2pt \hbox{$\mathrel{\scriptstyle\mathop{\kern0pt\prod}\limits_{#1}^{#2}}$}}

\font\eightrm=cmr8         \font\eighti=cmmi8
\font\eightsy=cmsy8        \font\eightbf=cmbx8
\font\eighttt=cmtt8        \font\eightit=cmti8
\font\eightsl=cmsl8        \font\sixrm=cmr6
\font\sixi=cmmi6           \font\sixsy=cmsy6
\font\sixbf=cmbx6\catcode`\@=11
\def\eightpoint{%
  \textfont0=\eightrm \scriptfont0=\sixrm \scriptscriptfont0=\fiverm
  \def\rm{\fam\z@\eightrm}%
  \textfont1=\eighti  \scriptfont1=\sixi  \scriptscriptfont1=\fivei
  \def\oldstyle{\fam\@ne\eighti}\let\old=\oldstyle
  \textfont2=\eightsy \scriptfont2=\sixsy \scriptscriptfont2=\fivesy
  \textfont\itfam=\eightit
  \def\it{\fam\itfam\eightit}%
  \textfont\slfam=\eightsl
  \def\sl{\fam\slfam\eightsl}%
  \textfont\bffam=\eightbf \scriptfont\bffam=\sixbf
  \scriptscriptfont\bffam=\fivebf
  \def\bf{\fam\bffam\eightbf}%
  \textfont\ttfam=\eighttt
  \def\tt{\fam\ttfam\eighttt}%
  \abovedisplayskip=9pt plus 3pt minus 9pt
  \belowdisplayskip=\abovedisplayskip
  \abovedisplayshortskip=0pt plus 3pt
  \belowdisplayshortskip=3pt plus 3pt 
  \smallskipamount=2pt plus 1pt minus 1pt
  \medskipamount=4pt plus 2pt minus 1pt
  \bigskipamount=9pt plus 3pt minus 3pt
  \normalbaselineskip=9pt
  \setbox\strutbox=\hbox{\vrule height7pt depth2pt width0pt}%
  \normalbaselines\rm}\catcode`\@=12

\font\san=cmssdc10
\def\ext{\hbox{\san \char3}}
\def\sym{\hbox{\san \char83}}
\def\pc#1{\tenrm#1\sevenrm}
\def\tx{\kern-1.5pt -}
\def\cqfd{\kern 2truemm\unskip\penalty 500\vrule height 4pt depth 0pt width
4pt\medbreak} 
\def\virg{\raise
.4ex\hbox{,}}
\def\ind{\par\hskip 1truecm\relax}
\def\End{\mathop{\rm End}\nolimits}
\def\im{\mathop{\rm Im}\nolimits}
\def\Ker{\mathop{\rm Ker}\nolimits}
\def\Pic{\mathop{\rm Pic}\nolimits}

\frenchspacing
\input amssym.def
\input amssym
\vsize = 25.6truecm
\hsize = 16truecm
%\hoffset = -.15truecm
\voffset = -.5truecm
\parindent=0cm
\baselineskip15pt
\overfullrule=0pt

\centerline{\bf Remark on a conjecture of Mukai}\smallskip
 \centerline{Arnaud {\pc BEAUVILLE}} 
\vskip1cm

{\bf Introduction}
\smallskip
\ind  The conjecture mentioned in the title appears actually as a question  in [M]
(Problem 4.11):
\par
{\bf  Conjecture}$.-$ {\it Let $C$ be a general curve, and  $E$ a
stable rank $2$ vector bundle on $C$ with $\det E=K_C$. The multiplication map
$\mu _E:\sym^2H^0(C,E)\rightarrow H^0(C,\sym^2E)$ is injective}.
  \ind Let ${\cal M}_K$ be the moduli space of stable rank $2$ vector
bundles $E$ on $C$ with $\det E=K_C$, and let ${\cal M}_K^n$ be the
subvariety of ${\cal M}_K$ parametrizing bundles with $h^0(E)= n$.
As explained in {\it loc. cit.}, the conjecture implies that ${\cal M}_K^n$ is smooth, of codimension ${1\over 2}n(n+1)$; this would give an analogue of the Brill-Noether theory for rank 2 vector bundles with canonical determinant.
\ind In this note we prove the conjecture in a very particular case:
\smallskip {\bf Proposition 1}$.-$ 
{\it The conjecture holds if
$h^0(E)\le 6$}.
\ind As a corollary we obtain that 
 ${\cal M}_K^n$ {\it is
smooth of codimension ${1\over 2}n(n+1)$ in ${\cal M}$ for} $n\leq 6$.
 Another consequence  is that Mukai's
conjecture holds for $g\le 9$ (see \S 2).
\ind The idea of the proof is the following. Put $n=h^0(E)$; assume for
simplicity that
$E$ is generated by its global sections, and that the map $\mu _E$
annihilates a non-degenerate (symmetric) tensor. This means that the image of the 
 map ${\bf P}_C(E)\rightarrow {\bf P}^{n-1}$ defined by the tautological
line bundle is contained in a smooth quadric; equivalently, the
map from $C$ into the Grassmann variety ${\bf G}(2,n)$ associated to
$E$ factors through the orthogonal Grassmannian ${\bf GO}(2,n)$. Now
for $n\le 6$, the restriction of the Pl\"ucker line bundle ${\cal O}_{\bf
G}(1)$ to  ${\bf GO}(2,n)$ is the sum of two (possibly equal) globally generated line
bundles; pulling back to $C$ we obtain a decomposition of $\det E=K_C$
which turns out to contradict Brill-Noether theory for a general curve.
\ind The last part of the argument breaks down for $n\ge 7$, since then the
restriction map $\Pic({\bf G}(2,n))\rightarrow \Pic({\bf GO}(2,n))$  is an
isomorphism. In fact we will show that the result cannot be improved without 
strengthening the hypotheses that we use on $C$ and $E$.
\vskip1truecm
{\bf 1. Proof of the main result}\smallskip 
\ind We will say that a curve $C$ is {\it Brill-Noether general} if for any line bundle
$L$ on $C$, the multiplication map
$$H^0(C,L)\otimes H^0(C,K_C\otimes L^{-1} )\longrightarrow H^0(C,K_C)$$is
injective. A general curve of given genus is Brill-Noether general; if $S$ is a K3
surface with $\Pic(S)={\bf Z}[H]$, a general
element of the linear system $|H|$ is Brill-Noether general [L].
\ind We will use this property in the following way:
 
{\bf Lemma 1}$.-$ {\it Let $C$ be a Brill-Noether general curve and $L,L'$ two
line bundles on $C$ such that
$h^0(K_C\otimes (L\otimes L')^{-1} )\ge 1$.
 Then :
\ind {\rm a)} The multiplication map $H^0(C,L)\otimes
H^0(C,L')\rightarrow  H^0(C,L\otimes L')$ is injective;
\ind {\rm b)} If $L'=L$, we have $h^0(L)\le 1$}.

{\it Proof} : Choose a non-zero section $s\in H^0(C,(K_C\otimes
(L\otimes L')^{-1} ))$, and consider the commutative diagram\vskip-10pt
$$\diagram{H^0(L)\otimes H^0(L')&\hfl{}{}& H^0(L\otimes L'))\cr
\vfl{1\otimes s}{}&&\vfl{}{ s}\cr
H^0(L)\otimes H^0(K_C\otimes L^{-1} )&\hfl{}{}& H^0(K_C)\ .\cr
}$$\vskip-10pt
By our hypothesis on $C$ the bottom horizontal map is injective; it
follows that the top one is injective.
\ind If $L'=L$, we get that the map $H^0(L)^{\otimes 2}\rightarrow
H^0(L^{\otimes 2})$ is injective; if $H^0(L)$ contains two linearly
independent elements $s,t$, the tensor $s\otimes t-t\otimes s$ is
non-zero and belongs to the kernel of that map, a contradiction.\cqfd
\smallskip 
\ind Let $E$ be a vector bundle on a curve $C$. 
If the multiplication map $\mu _E$ is injective, the same property holds for all
subbundles of $E$; thus $E$ satisfies

\ind  $(\star)$ {\it For every sub-line bundle $L\i E$, the  map $\mu
_L:\sym^2H^0(C,L)\rightarrow  H^0(C,L^2)$ is injective}.\smallskip 
\ind By [T], this condition is automatically satisfied if $C$ is general and  any
sub-line bundle of $E$ has degree $\le g+1$. Thus Proposition 1 is a consequence of
the following  more precise
result:
\smallskip 
{\bf Proposition 2}$.-$ {\it Let $C$ be a Brill-Noether general curve, and
$E$ a  rank $2$ vector bundle on $C$ with $\det E=K_C$, satisfying condition
$(\star)$. Then any non-zero tensor 
$\tau\in \sym^2H^0(C,E) $ such that $\mu _E(\tau )=0$ has rank} $>6$.

{\it Proof} : a) {\it The general set-up}
\ind  Let $C$ be a curve,  $E$ a rank 2 vector bundle on $C$
with $\det E=K_C$, $\tau $ an element of $\Ker \mu _E$ of rank $n$.
This means that we can find linearly independent elements $s_1,\ldots ,s_n$ of
$H^0(C,E)$ such that 
$\tau =s_1^2+\ldots +s_n^2$. Let $F\i E$ be the image of the map ${\cal
O}_C^n\rightarrow E$ defined by $s_1,\ldots ,s_n$. Then  $s_1,\ldots
,s_n$ are sections of $F$ which satisfy $s_1^2+\ldots +s_n^2=0$ in
$H^0(C,\sym^2F)$. By property $(\star)$  $F$ is a rank 2
subsheaf of $E$, with determinant $K_C(-A)$ for some effective divisor $A$.
\ind Let $P:={\bf P}_C(F)$, and let ${\cal O}_P(1)$ be the tautological line
bundle on  $P$. Through  the
 canonical isomorphisms 
$$H^0(P,{\cal O}_P(1))\iso H^0(C,F)\qquad H^0(P,{\cal O}_P(2))\iso
H^0(C,\sym^2F)$$ 
 the multiplication map $\mu _F$ is identified with $\mu
_{{\cal O}_P(1)}$; thus we can view $s_1,\ldots
,s_n$ as global sections of ${\cal O}_P(1)$, which generate ${\cal O}_P(1)$
and satisfy $s_1^2+\ldots +s_n^2=0$. In other words, 
 the image of the morphism $\varphi :P\rightarrow
{\bf P}^{n-1}$ defined by $(s_1,\ldots,s_n)$ is contained in the smooth
quadric $Q$ defined by
$X_1^2+\ldots +X_n^2=0$. 
\ind Let $\gamma $ be the map of $C$ into the Grassmann variety ${\bf
G}(2,n)$ associated to the surjective homomorphism ${\cal
O}_C^n\rightarrow F$; by definition $F$ is the pull back of the universal 
quotient bundle on 
${\bf G}(2,n)$. The determinant of that bundle is the Pl\"ucker line bundle
${\cal O}_{\bf G}(1)$ on ${\bf G}(2,n)$, so we get 
$\gamma ^*{\cal O}_{\bf G}(1)\cong K_C(-A)$.
\ind  For each $x\in C$,  the point $\gamma (x)\in{\bf
G}(2,n)$ corresponds to the line $\varphi ({\bf
P}(F_x))$ in ${\bf P}^{n-1}$. So the fact that the image of
$\varphi $ is contained in $Q$ means that $\gamma $ factors through 
the orthogonal grassmannian ${\bf GO}(2,n)$ of lines contained in $Q$:
$$\gamma : C\ \hfl{\gamma _{\bf O}}{}\ {\bf GO}(2,n)\mono {\bf
G}(2,n)\ .$$
 b) {\it The cases $n=4$ and $n=5$}
\ind Now we will assume $n\le 6$ and that $C$ is Brill-Noether general, and derive a
contradiction.   The case
$n\le 3$ is trivial (${\bf GO}(2,n)$ is empty). Consider the case $n=4$. 
 Then ${\bf GO}(2,4)$ parametrizes the lines in a smooth quadric
in ${\bf P}^3$; it has two components, both isomorphic to ${\bf P}^1$. Inside
${\bf G}(2,4)$, which is a quadric in ${\bf P}^5$, each of these components
is a conic. Thus $\gamma $ factors as
$$\gamma :C\ \hfl{\gamma _{\bf O}}{}\ {\bf P}^1\mono {\bf
G}(2,4)\ ,$$with ${\cal O}_{\bf G}(1)_{|{\bf P}^1}={\cal O}_{{\bf P}^1}(2)$. Thus
$L=\gamma _{\bf O}^*{\cal O}_{{\bf P}^1}(1)$ satisfies $L^2=K_C(-A)$ and \break
$h^0(L)\ge 2$, contradicting  Lemma 1 b).

\ind Suppose $n=5$. Let $V$ be a 4-dimensional vector space, with a non-degenerate
alternate form $\omega $ and an orientation $\ext^4V\iso{\bf C}$. The
orthogonal $W$ of $\omega $ in $\ext^2V$ is a 5-dimensional vector space
with a non-degenerate quadratic form given by the wedge product. The
isotropic 2-planes in $W$ are of the form $\ell \wedge\ell ^{\perp}$, for
$\ell\in {\bf P}(V) $. Thus  ${\bf GO}(2,5)\i {\bf G}(2,5)$ is
identified with
${\bf P}(V)$ embedded in ${\bf G}(2,W)$ by $\ell \mapsto \ell
\wedge\ell ^{\perp}$. The corresponding map from ${\bf P}(V)$ to ${\bf
P}(\ext^2W)$ is quadratic, so again ${\cal O}_{\bf G}(1)_{|{\bf P}(V)}={\cal O}_{{\bf
P}(V)}(2)$, and we conclude exactly as above.
\medskip 
c) {\it The case $n=6$}
\ind  Let again $V$ be a 4-dimensional vector space with an
orientation. Then $\ext^2V$ is a 6-dimensional quadratic  vector space; the
corresponding quadric in ${\bf P}(V)$ is the Grassmannian ${\bf G}(2,V)$.
The lines contained in ${\bf G}(2,V)\i{\bf P}(\ext^2V)$ correspond to the 2-planes
$\ell
\wedge H\i\ext^2V$ where
$\ell\i H\i V
$,
$\dim\ell =1$, $\dim H=3$; thus ${\bf GO}(2,6)\i {\bf G}(2,6)$ is
identified with the incidence variety 
$Z \i{\bf P}(V)\times {\bf P}(V^*)$,
 embedded in ${\bf G}(2,V)$ by $(\ell ,H)\mapsto
\ell\wedge H$. Fixing $\ell $ or $H$ one sees that the restriction of ${\cal
O}_{\bf G}(1)$ to  $Z$ is the pull back of ${\cal O}_{\bf P}(1)\boxtimes
{\cal O}_{\bf P}(1)$.
\ind The map $\gamma _{\bf O}:C\rightarrow Z$ gives by projection two maps
$u:C\rightarrow {\bf P}(V)$ and $u':C\rightarrow {\bf P}(V^*)$.  Put $L=u^*{\cal
O}_{{\bf P}(V)}(1)$, $L'=u'^*{\cal O}_{{\bf P}(V^*)}(1)$. Then  $ L\otimes L' \cong
K_C(-A)$, and we deduce from $u$ and $u'$ two homomorphisms
$$v:V^*\rightarrow H^0(C,L) \qquad v':V\rightarrow H^0(C,L')\ .$$
\ind Consider the map
$$V\otimes V^*\ \hfl{v\otimes v'}{}\ H^0(L)\otimes
H^0(L')\longrightarrow H^0(L\otimes L')\ ;$$ the fact that 
$(u,u'):C\rightarrow {\bf P}(V)\times {\bf P}(V^*)$ factors through $Z$ means that
the identity element of
$V\otimes V^*\cong \End(V)$ goes to zero in $H^0(L\otimes L')$, and
therefore already in $H^0(L)\otimes H^0(L')$ by Lemma 1 a). 
\ind Put $K=\Ker u$, $K'=\Ker u'$. The kernel of $v\otimes v'$ is $K\otimes
V^*+V\otimes K'$; any element of this kernel has rank $\leq \dim K+\dim K'$. Since
the identity tensor has rank 4, we get $\dim K+\dim K'\geq 4$.
\ind Suppose $\dim K=\dim K'=2$. Identifying again $V\otimes V^*$ to $ \End(V)$,
we can write ${\rm Id}_V=p+q$, with $\im p\i K$ and $\Ker q\supset K'^{\perp}$.
This implies that $p$ and $q$ are orthogonal projectors, and therefore that
$K'=K^{\perp}$.   Then ${\bf P}(K)\times {\bf P}(K^{\perp})$ is contained in $Z$,
and $\gamma _{\bf O}$ factors as
$$C\ \hfl{(u,u')}{}\ {\bf P}(K)\times {\bf P}(K^{\perp})\mono {\bf GO}(2,6)\ .$$
 Let $(\ell,H)\in {\bf P}(K)\times {\bf P}(K^{\perp})$; then $\ell\i K$ and $H\supset
K$, so  the line ${\bf P}(\ell\wedge H)$ in ${\bf P}(\ext^2V)={\bf P}^5$ contains the
point
$p:={\bf P}(\ext^2K)$. In other words, the lines in $Q$ parametrized by $C$ all pass
through
$p$, hence are contained in the tangent hyperplane $T_p(Q)$. Therefore the scroll
$\varphi (P)$ in
${\bf P}^5$ is contained in $T_p(Q)$; this is impossible because $\varphi $ is defined
by 6 linearly independent sections of ${\cal O}_P(1)$.
\ind Suppose now $\dim K=3$, so that the image of $C$ in ${\bf P}(V^*)$ is the
point ${\bf P}(K^{\perp})$. Then $\gamma_{\bf O} $ factors as
$$C\qfl{u} {\bf P}(K)\mono {\bf GO}(2,6)\ ;$$the image of  ${\bf P}(K)$ in ${\bf
GO}(2,6)$ is the family of lines contained in ${\bf
P}(\ext^2K)\i {\bf P}(\ext^2V)$ $={\bf P}^5$. Thus $\varphi $ maps $P$ to the
projective plane ${\bf P}(\ext^2K)$, which is again impossible. The same argument
applies when $\dim K'=3$.\cqfd

\vskip1truecm
{\bf 2. Consequences and comments}\smallskip 
\ind The following lemma must be well-known:\par
{\bf Lemma 2}$.-$ {\it Let $C$ be a Brill-Noether general curve, and $E$  a
semi-stable rank $2$ vector bundle
 on $C$ with $\det E=K_C$. Then  $h^0(E)\le {g\over 2}+2$}.\par
{\it Proof} : Let $L$ be a sub-line bundle of $E$ of maximal degree. The exact
sequence
$$0\rightarrow L\rightarrow E\rightarrow K_C\otimes
L^{-1}
\rightarrow 0$$gives
$h^0(E)\le h^0(L)+h^0(K_C\otimes
L^{-1})$. Since 
$h^0(L)\,h^0(K_C\otimes L^{-1})\le g$ and $h^0(L)\le $ $h^0(K_C\otimes L^{-1})$
by semi-stability, the required inequality holds unless  $h^0(L)$  is $0$ or $1$. In
that case we have $$h^0(L)+h^0(K_C\otimes L^{-1})=2h^0(L)+g-1-\deg L\le g+1-\deg
L$$and $\deg L\ge {g\over 2}-1$ by [N], hence 
$h^0(E)\le  {g\over 2}+2$.\cqfd
\ind In particular $g\le 9$ guarantees $h^0(E)\le 6$, hence:\par
{\bf Corollary}$.-$ {\it Mukai's conjecture holds for $g\le 9$}.
\smallskip 
\ind We will show that Proposition 2 cannot be improved without further hypotheses:
\par
{\bf Proposition 3}$.-$ {\it For each integer $s\ge 5$, there exists a Brill-Noether
general curve $C$ of genus $g=2s$ and a stable rank $2$ vector bundle $E$ on $C$
with
$\det E=K_C$, satisfying $(\star)$, such that $\Ker\mu _E$ contains a tensor of
rank $7$}.\par\smallskip 
{\it Proof} : We choose a K3 surface $S$ with $\Pic(S)={\bf Z}[H]$, $(H^2)=4s-2$, and
take for
$C$ a general element of the linear system $|H|$. Then $C$ has genus $2s$ and is
Brill-Noether general;  in particular it admits a finite set ${\cal P}$ of line
bundles $L$ with
$h^0(L)=2$, $\deg L=s+1 $. Choose one of these line bundles, say $L$; the {\it
Lazarsfeld bundle} $E_L$ is defined by the exact sequence 
$$0\rightarrow E_L^*\longrightarrow H^0(C,L)\otimes {\cal O}_S\longrightarrow
L\rightarrow 0\ .$$
 It turns out that the restriction $E$ of $E_L$ to $C$ does not depend on the choice of
$L\in {\cal P}$; we have $\det E=K_C$, $h^0(E)=s+2$ (see [V]), and for each
$L\in{\cal P}$ an exact sequence
$$0\rightarrow L\rightarrow E\rightarrow K_C\otimes
L^{-1}\rightarrow 0\ .$$

{\bf Lemma 3}$.-$ {\it Let $M$ be a sub-line bundle of $E$. Then either $M\in{\cal
P}$, or $h^0(M)\le 1$ and $\deg M\le s$}.\par
{\it Proof} : Put $M':=K_C\otimes M^{-1}$; we may assume that the  $E/M$ is
 isomorphic to $M'$, so that
$$s+2=h^0(E)\le
h^0(M)+h^0(M')=2h^0(M)+g-1-\deg M=2h^0(M')+g-1-\deg M'\ .$$
As in the proof of Lemma 2 this implies that either $h^0(M)$ or $h^0(M')$ is $\le
2$. 
\ind If $h^0(M)\le 2$, we get  $\deg M\le s+1$, hence
either $M\in {\cal P}$, or $\deg M\le s$ and $h^0(M)\le 1$. 
\ind Assume $h^0(M')\le 2$, so that $h^0(M)\ge s$. We get  again $\deg M'\le s+1$,
hence $\deg M\ge 3s-3$.
On the
other hand the exact sequence preceding the Lemma shows that $M$ injects into
$K_C\otimes L^{-1}$ for each $L\in{\cal P}$; since ${\rm Card}({\cal P})\ge 2$ the
inclusion must be strict, so $\deg M<\deg K_C\otimes
L^{-1}=3s-3$, a contradiction.\cqfd

\ind It follows that $E$ is stable and satisfies condition $(\star)$. We observe that
$\sym^2E$ is isomorphic to 
$K_C\otimes {\cal E}nd_0(E)$, where ${\cal E}nd_0(E)$is the sheaf  of trace 0
endomorphisms of $E$; thus  Serre duality and the stability of $E$ imply
$h^1(\sym^2E)=h^0({\cal E}nd_0(E))=0$, hence $h^1(\sym^2E)=3g-3$ by
Riemann-Roch.
\ind In
$\sym^2H^0(E)$ the locus of tensors of rank $\le 7$ has dimension $7(s-1)$, while
the kernel of $\mu _E$ has codimension $\le 6s-3$. Therefore the intersection of
these subvarieties is not reduced to $0$ -- in fact it has dimension $\ge s-4$. By
Proposition 2 all the tensors in this intersection have rank $7$.\cqfd

\vskip2truecm
\centerline{ REFERENCES} \vglue15pt\baselineskip12.8pt
\def\num#1{\smallskip\item{\hbox to\parindent{\enskip [#1]\hfill}}}
\parindent=1.3cm 
\num{L} R. {\pc LAZARSFELD}: {\sl Brill-Noether-Petri without degenerations}. J.
Differential Geom. {\bf 23} (1986), no. 3, 299--307.
\num{M} S. {\pc MUKAI}: {\sl Vector bundles and Brill-Noether theory}. Current topics
in complex algebraic geometry,  145--158; Math. Sci. Res. Inst. Publ. {\bf 28},
Cambridge Univ. Press,  1995. 

\num{N} M. {\pc NAGATA}:
{\sl On self-intersection number of a section on a ruled surface}.
Nagoya Math. J. {\bf 37} (1970), 191--196.
\num{T} M. {\pc TEIXIDOR}: {\sl Injectivity of the symmetric map for line bundles}.
Manu\-scripta Math. {\bf 112} (2003), no. 4, 511--517. 
\num{V} C. {\pc VOISIN}: {\sl Sur l'application de Wahl des courbes satisfaisant la
condition de Brill-Noether-Petri}. Acta Math. {\bf  168} (1992), no. 3-4, 249--272. 
\vskip1cm
\def\pc#1{\eightrm#1\sixrm}
\hfill\vtop{\eightrm\hbox to 5cm{\hfill Arnaud {\pc BEAUVILLE}\hfill}
 \hbox to 5cm{\hfill Institut Universitaire de France\hfill}\vskip-2pt
\hbox to 5cm{\hfill \&\hfill}\vskip-2pt
 \hbox to 5cm{\hfill Laboratoire J.-A. Dieudonn\'e\hfill}
 \hbox to 5cm{\sixrm\hfill UMR 6621 du CNRS\hfill}
\hbox to 5cm{\hfill {\pc UNIVERSIT\'E DE}  {\pc NICE}\hfill}
\hbox to 5cm{\hfill  Parc Valrose\hfill}
\hbox to 5cm{\hfill F-06108 {\pc NICE} Cedex 02\hfill}}
\end